\documentclass[12pt]{amsart}
\usepackage{amsmath,amsthm,amssymb, amscd, amsfonts}

\newcommand{\uc}{{\mathcal U}}
\newcommand{\acw}{\widetilde{\mathcal A}}
\newcommand{\bc}{{\mathcal B}}
\newcommand{\Ss}{{\mathcal S}}
\newcommand{\otravez}{{1\le i \le \theta}}
\newcommand{\otraavez}{{1\le i, j \le \theta}}
\newcommand{\otrvez}{{1\le i \neq  j \le \theta}}

\newcommand{\ttheta}{\widehat{\theta}}

\newcommand{\ku}{{\Bbbk}}
\newcommand{\Z}{{\mathbb Z}}
\newcommand{\N}{{\mathbb N}}
\newcommand{\Q}{{\mathbb Q}}
\newcommand{\kv}{{\Q (v)}}

\newcommand{\C}{{\mathbb C}}
\newcommand{\ac}{{\mathcal A}}

\newcommand{\VGamma}{\widehat{\Gamma}}

\newcommand{\ydg}{{}_{\ku \Gamma}^{\ku \Gamma}\mathcal{YD}}
\newcommand{\ydh}{{}_{G(A)}^{G(A)}\mathcal{YD}}
\newcommand{\ydhh}{{}_H^{H}\mathcal{YD}}
\newcommand{\toba}{{\mathfrak B}}
\newcommand{\wtoba}{\widetilde {\mathfrak B}}

\newcommand{\tobak}{\toba(W)_{\ku}}
\newcommand{\wtobak}{\wtoba(W)_{\ku}}
\newcommand{\wtobaka}{\wtoba(W)_{\ac}}

\newcommand{\wl}{\widetilde L}
\newcommand{\wll}{\widetilde \ell}

\numberwithin{equation}{section}
\newcommand{\height}{\mbox{\rm \, ht\, }}

\newcommand{\End}{\mbox{\rm End\,}}

\newcommand{\id}{\mathop{\rm id\,}}

\newcommand{\ord}{\mathop{\rm ord}}

\newcommand{\ad}{\mbox{\rm ad\,}}

\newcommand{\gr}{\mbox{\rm gr\,}}

\theoremstyle{plain}

\newtheorem{teo}{Theorem}[section]
\newtheorem{lema}[teo]{Lemma}
\newtheorem{cor}[teo]{Corollary}
\newtheorem{prop}[teo]{Proposition}

\theoremstyle{definition}
\newtheorem{defi}[teo]{Definition}

\newtheorem{obs}[teo]{Remark}

\def\pf{\begin{proof}}
\def\epf{\end{proof}}

\begin{document}

To appear in Crelle Journal. \bigbreak

\renewcommand{\baselinestretch}{1.2}
\renewcommand{\thefootnote}{}
\thispagestyle{empty}
\title{A characterization
of quantum groups}
\author{ Nicol\'as Andruskiewitsch and Hans-J\"urgen Schneider}
\address{Facultad de Matem\'atica, Astronom\'\i a y F\'\i sica\\
Universidad Nacional de C\'ordoba \\ (5000) Ciudad Universitaria
\\C\'ordoba \\Argentina}
\email{andrus@mate.uncor.edu}
\dedicatory{In  memoriam Peter Slodowy}
\begin{abstract}
We classify pointed Hopf algebras with finite Gelfand-Kirillov
dimension, which are domains, whose groups of group-like elements are
finitely generated and abelian,
and whose infinitesimal braidings are positive.
\end{abstract}
\thanks{
 \\
This work was partially supported by ANPCyT, Agencia C\'ordoba Ciencia,
CONICET, the Graduiertenkolleg of the Math. Institut
(Universit\"at \ M\"unchen)  and Secyt (UNC)}
\address{Mathematisches Institut\\
Universit\"at \ M\"unchen\\
Theresienstra\ss e 39\\
D-80333 M\"unchen\\
Germany}
\email{hanssch@rz.mathematik.uni-muenchen.de}
\maketitle
\section*{Introduction}
Since the appearance of quantum groups \cite{KR, Sk, Dr, Ji}, there were
many attempts to define them  intrinsically.
Important descriptions of the so-called "nilpotent" parts were given by
Ringel \cite{Ri}, Lusztig \cite{L3} and Rosso \cite{Ro3}.
However, the question of finding an abstract characterization
of the  quantized enveloping algebras remained open.

In the main Theorem \ref{fingrowth-lifting} of this paper, we classify all Hopf algebras
over an algebraically closed field of characteristic 0 which are
\begin{itemize}
\item pointed, that is all their simple comodules are one-dimensional, and have a finitely generated
abelian group of group-like elements,
\item domains of finite Gelfand-Kirillov dimension, and
\item have positive infinitesimal braiding (see Section \ref{corfilt}).
\end{itemize}
The first  two conditions are natural.
The positivity condition should
be related to the existence of a real involution.

In Theorem \ref{construction}, we describe these Hopf algebras by generators and
relations. They are natural generalizations of quantized enveloping algebras with
positive parameter.
To prove our main Theorem,
we combine the lifting method for pointed Hopf algebras
\cite{AS1, AS2, AS4} with a characterization obtained by Rosso
of the "nilpotent part" of a quantized enveloping algebra
in terms of  finiteness of the Gelfand-Kirillov dimension  \cite{Ro3}.

\medbreak

Among the main differences between the new Hopf algebras
and multi-parametric quantized enveloping algebras,
let us mention that
we have one parameter of deformation
for each connected component of the Dynkin diagram
(this is  explained as follows:
in the "classical limit", one may
have different scalar multiples of
the Sklyanin brackets in the different connected components)
and  linking relations, see \eqref{relations3},
generalizing the classical relations $E_i F_j - F_j E_i
= \delta_{ij} (K_i - K_i^{-1})$.

\medbreak
Note that we are not assuming that the Hopf algebras
have an {\em a-priori} assigned
Dynkin diagram as in \cite{W1, KW};
it comes from our hypothesis,
and here is where we rely on Rosso's result \cite{Ro3}.

\medbreak
The article is organized as follows. Section 1 contains
preliminary material. In Section 2, we collect some well-known
facts about quantum groups, and give details of some proofs
when they are not easily available in the literature.
Section 3 contains a technical result on the coradical
filtration of certain Hopf algebras, generalizing an idea of
Takeuchi. In Section 4, we construct a new family
of pointed Hopf algebras with generic braiding
and establish the main basic properties of them.
The approach is similar to \cite{AS4} but instead
of dimension arguments, we use the technical results
on the coradical filtration obtained in Section 3;
these results should be useful also for other classes of
Hopf algebras. In Section 5, we prove our main Theorem.
A key point is Lemma \ref{fingrowth-nichols},
which implies that a wide class of Hopf algebras
with finite Gelfand-Kirillov dimension is generated by
group-like and skew-primitive elements.

\subsubsection*{Acknowledgments} We thank N. Reshetikhin for reviving our interest in
this
question and J. Alev for interesting conversations.
This paper was begun during visits to
the MSRI, in the framework of the full-year
Program on Non-commutative Algebra (August 1999
- May 2000). We thank the
organizers for the kind invitation and the MSRI for the excellent working conditions.
Part of the work of the first author was done during a visit
to the University of Rheims (October 2001 - January 2002);
he is very grateful to J. Alev for his kind hospitality.

\section{Preliminaries}\label{corfilt}

\subsection*{Notation}
Let $\ku$ be an algebraically closed field of characteristic 0.
Our references  are \cite{M}, \cite{Sw} for Hopf algebras;
\cite{KL}, for growth of algebras and Gelfand-Kirillov
dimension; and \cite{AS5} for pointed Hopf algebras.
We use standard notation for Hopf algebras:
$\Delta$, $\Ss$, $\epsilon$, denote respectively the
comultiplication, the antipode, the counit; we use
a short version of Sweedler's
notation: $\Delta (x) = x_{(1)} \otimes x_{(2)}$, $x$
in a coalgebra $C$.

\medbreak Let $H$ be a Hopf algebra with bijective antipode. We
denote by $G(H)$  the group of group-like elements of $H$; and by
${\mathcal P}_{g,h}(H)$ the space of $g$, $h$ skew-primitive
elements $x$ of $H$, that is with $\Delta(x) = g \otimes x + x
\otimes h$, where $g, h\in G(H)$; then ${\mathcal P}(H) :=
{\mathcal P}_{1,1}(H)$. The braided category  of Yetter-Drinfeld
modules over  $H$ is denoted by $\ydhh$, {\it cf.} the conventions
of \cite{AS5}.

\medbreak The  adjoint representation $\ad$ of a Hopf algebra $A$
on itself is given by $\ad x(y)= x_{(1)} y \Ss(x_{(2)})$. If $R$
is a braided Hopf algebra in $_{H}^{H}\mathcal{YD}$ then there is
a braided adjoint representation $\ad_c$ of $R$ on itself given by
$\ad_c x(y)=\mu(\mu\otimes\Ss)(\id\otimes
c)(\Delta\otimes\id)(x\otimes y),$ where $\mu$ is the
multiplication and $c\in\End(R\otimes R)$ is the braiding. If
$x\in {\mathcal P} (R)$ then the braided adjoint representation of
$x$ is $ \ad_cx(y)=\mu(\id-c)(x\otimes y) =: [x, y]_{c}$. The
element $[x, y]_{c}$ defined by the second equality for any $x$
and $y$, regardless of whether $x$ is primitive, will be called a
braided commutator. When $A = R\# H$ (where $\#$ stands for
Radford biproduct or bosonization), then for all $b,d\in R$,
$\ad_{(b\# 1)}(d\# 1)=(\ad_cb(d))\# 1$.

\medbreak
If $\Gamma$ is an abelian group, we denote by $\widehat{\Gamma}$ the
group of characters of $\Gamma$.
If $V$ is a $\ku \Gamma$-module (resp., $\ku \Gamma$-comodule),
then we denote $V^{\chi} :=
\{v\in V: h.v = \chi(h)v, \forall h\in \Gamma\}$, $\chi \in \widehat{\Gamma}$;
resp., $V_{g} := \{v\in V: \delta(v) = g\otimes v\}$, $g\in G$.
A Yetter-Drinfeld module over $\Gamma$ is
$\ku \Gamma$-module $V$ which is also a $\ku \Gamma$-comodule, and  such that
each homogeneous component $V_{g}$, $g\in \Gamma$, is a $\ku\Gamma$-submodule.
Thus, a vector space $V$ provided with a direct sum
decomposition $V = \oplus_{g\in G, \chi\in \widehat{\Gamma}}
V^{\chi}_{g}$
is a Yetter-Drinfeld module over $H = \ku\Gamma$.

\subsection*{Braided vector spaces}
A  braided vector space
$(V, c)$ is a finite-dimensional vector space provided with an isomorphism
$c: V\otimes V\to V\otimes V$ which is a solution of the braid equation,
that is $(c\otimes \id) (\id\otimes c) (c\otimes \id) =
(\id\otimes c) (c\otimes \id) (\id\otimes c).$
Examples of braided vector spaces are
Yetter-Drinfeld modules: if $V\in \ydhh$, then
$c: V\otimes V \to V\otimes V$, $c(v \otimes w) = v_{(-1)}.w \otimes v_{(0)}$, is a
solution of the braid equation.

\begin{defi}\label{defis}
Let $(V, c)$ be a finite-dimensional braided vector space.
We shall say that the braiding $c: V\otimes V\to
V\otimes V$ is {\it diagonal} if there exists a basis
$x_{1}, \dots, x_{\theta}$ of $V$ and non-zero scalars $q_{ij}$
such that $c(x_{i}\otimes x_{j}) = q_{ij} x_{j}\otimes x_{i}$,
$\otraavez$. The matrix $(q_{ij})$
is called the matrix of the braiding.

\medbreak
Furthermore, we shall say that a diagonal braiding with matrix $(q_{ij})$
is {\it indecomposable}  if for all $i \neq j$,
there exists a sequence $i= i_{1}$, $i_{2}$, \dots, $i_{t} = j$
of elements of $\{1, \dots, \theta\}$ such that
$q_{i_s, i_{s+1}}q_{i_{s+1}, i_{s}} \neq 1$, $1\le s \le t-1$.
Otherwise, we say that the matrix is decomposable.

\medbreak
We  attach a graph to a diagonal braiding in the following way.
The vertices of the graph are the elements of $\{1, \dots, \theta\}$, and
there is an edge between $i$ and $j$ if they are different and $q_{ij}q_{ji} \neq 1$.
Thus, ``indecomposable" means that the corresponding graph is connected.
The components of the matrix are the principal
submatrices corresponding to the connected components of the graph.
If $i$ and $j$ are vertices in the same connected component, then we write $i\sim j$.
We shall denote by $\mathcal  X$ the set of connected components
of the matrix $(q_{ij})$. If $I\in \mathcal  X$, then $V_I$
denotes the subspace of $V$ spanned by $x_i$, $i\in I$.

\medbreak
We shall say that a   braiding $c$ is {\it generic}
if it is diagonal  with matrix $(q_{ij})$
where $q_{ii}$ is not a root of 1, for any $i$.

\medbreak Let $\ku = \C$.
We shall say that a   braiding $c$ is {\it positive}
if it is generic  with matrix $(q_{ij})$
where $q_{ii}$ is a positive real number, for all $i$.

\medbreak
We shall say that a  diagonal braiding $c$  with matrix $(q_{ij})$
is of {\it Cartan type} if  $q_{ii} \neq 1$ for all $i$, and
there are integers $a_{ij}$ with $a_{ii} = 2$,
$1\le i \le \theta$, and $0 \le -a_{ij} < \ord q_{ii}$ (which could be infinite),
$1\le i \neq j \le \theta$, such that
$q_{ij}q_{ji} = q_{ii}^{a_{ij}}$ for all $i$ and $j$.
Since clearly
$a_{ij} =0$ implies that $a_{ji} =0$ for all $i\neq j$,
$(a_{ij})$ is a generalized Cartan matrix. This generalizes the definition in
\cite[p. 4]{AS2}. In this case, the braiding is indecomposable if and only if the
corresponding Cartan matrix is indecomposable.
We shall also denote by $\mathcal  X$ the set of connected components
of the Dynkin diagram corresponding to the matrix $(a_{ij})$; clearly,
this agrees with the previous convention.

\medbreak
Let $(V,c)$ be a braided vector space of Cartan type with
generalized Cartan matrix $(a_{ij})$.
We say that $(V,c)$ is of {\it DJ-type} (or Drinfeld-Jimbo type)
if there exist positive integers $d_{1}, \dots, d_{\theta}$ such that
\begin{flalign}
 &\text{For all } i,j,\  d_{i} a_{ij} = d_{j} a_{ji} \text{ (thus } (a_{ij})
\text{ is symmetrizable).}\\
&\text{For all } I \in \mathcal  X, \text{there exists  } q_I \in
\ku \text{, which is not a root of unity,  such that
}\end{flalign} $$q_{ij} = q_I^{d_{i}a_{ij}} \text{ for all } i \in
I, 1\le j \le \theta.$$
 In particular, $q_{ij} = 1$ if $i\not \sim j$.
\end{defi}

\begin{lema}\label{uniq-basis}
Let $(V, c)$ be a finite-dimensional braided vector space
with diagonal braiding and matrix $(q_{ij})$, with respect to a basis
$x_{1}, \dots, x_{\theta}$ of $V$.
If there exist another basis
$y_{1}, \dots, y_{\theta}$ of $V$ and non-zero scalars $p_{ij}$
such that $c(y_{i}\otimes y_{j}) = p_{ij} y_{j}\otimes y_{i}$,
$\otraavez$, then there exists $\sigma \in \mathbb S_{\theta}$
such that $q_{ij} = p_{\sigma(i) \sigma(j)}$, $\otraavez$.
\end{lema}

\pf Let $(\alpha_{hr})$ be the transition matrix:
$y_r = \sum_{1 \le h \le \theta} \alpha_{hr} x_h$.
Then
$$
\sum_{1 \le h, l \le \theta}  p_{rs} \alpha_{hr}  \alpha_{ls} x_l \otimes x_h =
c(y_r \otimes y_s) = \sum_{1 \le h, l \le \theta}  \alpha_{hr}  \alpha_{ls}
\, c (x_h \otimes x_l) =
\sum_{1 \le h, l \le \theta} q_{hl} \alpha_{hr}  \alpha_{ls} x_l \otimes x_h.
$$
Hence $ p_{rs} \alpha_{hr}  \alpha_{ls} =  q_{hl} \alpha_{hr}  \alpha_{ls}$,
for all $1\le h, l, r, s \le \theta$.
Since the transition matrix is invertible,
there exists  $\sigma \in \mathbb S_{\theta}$ such that
$\alpha_{h \sigma(h)} \neq 0$, for all $1\le h \le \theta$.
The Lemma follows. \epf

\begin{lema}\label{det-basis}
Let $(V, c)$ be a finite-dimensional braided vector space
with generic braiding of Cartan type and matrix $(q_{ij})$,
with respect to a basis  $y_{1}, \dots, y_{\theta}$.
Assume that the braiding $c$ arises from a Yetter-Drinfeld
module structure on $V$ over an abelian group $\Gamma$.
Then there exist
$g_{1}, \dots, g_{\theta} \in
\Gamma$, $\chi_{1}, \dots, \chi_{\theta} \in
\widehat{\Gamma}$ and a basis $x_{1}, \dots, x_{\theta}$
such that
$x_i \in V_{g_{i}}^{\chi_{i}}$
and $c(x_i \otimes x_j) = q_{ij} x_j \otimes x_i$, $\otraavez$.
\end{lema}

\pf There exists a basis
$x_{1}, \dots, x_{\theta}$ of $V$ and $g_1, \dots, g_{\theta} \in \Gamma$
such that $x_i \in V_{g_i}$, $\otravez$.
Let $(\alpha_{hr})$ be the transition matrix:
$y_r = \sum_{1 \le h \le \theta} \alpha_{hr} x_h$.
Then
$$
\sum_{1 \le h \le \theta}  q_{rs} \alpha_{hr}  y_s \otimes x_h =
c(y_r \otimes y_s) = \sum_{1 \le h \le \theta}  \alpha_{hr}
\, c (x_h \otimes y_s) =
\sum_{1 \le h \le \theta} \alpha_{hr} \, g_h \cdot y_s \otimes x_h.
$$
Hence $ q_{rs} \alpha_{hr} y_s =  \alpha_{hr} \, g_h \cdot y_s$,
for all $1\le h, r, s \le \theta$; this implies that
the subgroup $\Gamma_0$ of $\Gamma$ generated
by  $g_1, \dots, g_{\theta}$ acts diagonally on $V$.
We can then refine the choice of the basis
$x_{1}, \dots, x_{\theta}$ and assume that
$x_i \in V_{g_i}^{\chi_i}$ for some  $\chi_1, \dots, \chi_{\theta} \in
\widehat{\Gamma_0}$;
by Lemma \ref{uniq-basis} we can assume (up to a permutation)
that
$c(x_i \otimes x_j) = q_{ij} x_j \otimes x_i$, $\otraavez$.
We claim that $(g_i, \chi_i) = (g_j, \chi_j)$ implies $i = j$.
If not, consider the subspace $W$ spanned by $x_i$ and $x_j$;
note that $q_{ij} = q_{ji} = q_{ii}$, hence
$q_{ii}^2 = q_{ij}  q_{ji} = q_{ii}^{a_{ij}}$.
Since the braiding is generic, $a_{ij} = 2$, a contradiction.
This proves the claim. Since the isotypic component
$V^{\chi}_g$ is $\Gamma$-stable, for any $g\in \Gamma$
and $\chi \in \widehat{\Gamma_0}$,
and  $\dim V^{\chi}_g \le 1$ by the claim, we see that
$\Gamma$ acts diagonally on $V$.
\epf

\subsection*{Nichols algebras}
Let $V\in \ydhh$. A braided graded Hopf algebra $R=
\oplus_{n\ge 0} R(n)$ in $\ydhh$
is called a {\it Nichols algebra} of $V$ if
$\ku\simeq R(0)$ and  $V\simeq R(1)$ in $\ydhh$,
and
\begin{enumerate}
\item $P(R) = R(1)$,
\item $R$ is generated as an algebra by $R(1)$.
\end{enumerate}

The Nichols algebra of $V$ exists and is unique up to isomorphisms;
it will be denoted by $\toba(V)$. It
depends, as an  algebra and coalgebra, only on the underlying braided vector space $(V,
c)$.
The underlying algebra is called a  {\it quantum symmetric algebra} in \cite{Ro3}.
We shall identify  $V$ with the subspace of homogeneous
elements of degree one in $\toba(V)$. See \cite{AS5} for more details
and some historical references.

\medbreak
Given a braided vector space of any of the types in Definition \ref{defis},
we will say that its Nichols algebra is of the same type.

\begin{lema}\label{conncompo} \cite[Lemma 4.2]{AS2}.
Let $V$ be a finite-dimensional Yetter-Drinfeld module over an abelian group. Let
$\mathcal  X
=\{I_1, \dots, I_N\}$ be a numeration of the set of connected components.
Then $\toba(V) \simeq \toba(V_{I_1}) \,\underline{\otimes} \dots \underline{\otimes} \,
\toba(V_{I_N})$
as braided Hopf algebras with the braided tensor product algebra structure
$\underline{\otimes}$ .    \qed \end{lema}

\subsection*{Lifting method for pointed Hopf algebras}
Recall that a Hopf algebra $A$ is pointed if
any irreducible $A$-comodule is one-dimensional. That is, if the coradical $A_{0}$
equals the group algebra $\ku G(A)$.

\medbreak Let $A$ be a pointed Hopf algebra let $A_{0} = \ku G(A)
\subseteq A_{1} \subseteq \dots$ be the coradical filtration and
let $\gr A = \oplus_{n\ge 0} \gr A(n)$ be the associated graded
coalgebra, which is  a graded Hopf algebra \cite{M}. The graded
projection $\pi: \gr A \to \gr A(0) \simeq \ku G(A)$ is a Hopf
algebra map and a retraction of the inclusion.
 Let $R = \{a\in A: (\id\otimes \pi)\Delta (a) = a\otimes 1\}$
 be the algebra of coinvariants of $\pi$; $R$ is a  braided Hopf algebra in
the category $\ydh$ of Yetter-Drinfeld modules
over $\ku G(A)$  and $\gr A$ can be reconstructed  from
$R$ and $\ku G(A)$ as a bosonization: $ \gr A \simeq R \# \ku G(A)$. Moreover,
$R = \oplus_{n\ge 0} R(n)$, where $R(n) = \gr A(n) \cap R$
is a graded braided Hopf algebra. We then have several invariants of our initial
pointed Hopf algebra $A$:
The graded braided Hopf algebra $R$; it is called the {\it diagram} of $A$.
 The braided vector space $(V, c)$, where  $V := R(1) = P(R)$
and $c: V\otimes V\to V\otimes V$ is the braiding in $\ydh$. It will be called
the {\it infinitesimal braiding} of $A$.
 The dimension of $V = P(R)$,  called the {\it rank} of $A$, or of $R$.
The subalgebra $R'$ of $R$ generated by $R(1)$, which  is the
Nichols algebra of $V$: $R' \simeq \toba(V)$.
See \cite{AS5} for more details.

\section{Nichols algebras of Cartan type}\label{nichols-cartantype}

\subsection*{Nichols algebras of diagonal type}

In this section, $(V, c)$ denotes a finite-dimensional braided vector
space; we assume that the braiding $c$ is diagonal with matrix $(q_{ij})$,
with respect to a basis $x_{1}, \dots x_{\theta}$.

Let $\Gamma$ be the free abelian group of rank $\theta$
with basis $g_{1}, \dots, g_{\theta}$.
We define characters $\chi_{1}, \dots, \chi_{\theta}$ of $\Gamma$ by
$$
\chi_{i} (g_{j}) = q_{ji}, \qquad 1\le i,j\le \theta.
$$
We consider $V$ as a Yetter-Drinfeld module over $\ku \Gamma$ by defining
$x_{i} \in V^{\chi_{i}}_{g_{i}}$, for all $i$.

\medbreak
We begin with some relations that hold in any Nichols algebra
of diagonal type.

\begin{lema}\label{primitivos}  Let $R$ be a braided Hopf algebra in $\ydg$, such that
$V \hookrightarrow P(R)$.

(a). If $q_{ii}$ is a root of 1 of order $N > 1$ for some $i\in \{1,
\dots, \theta\}$, then $x_{i}^N \in P(R)$.

(b). Let $i\neq j \in \{1,\dots, \theta\}$ such that $q_{ij}q_{ji} = q_{ii}^{1 -r}$,
where $r$ is an integer such that $ 0\le r-1 < \ord q_{ii}$
(which could be infinite). Then $(\ad x_{i})^{r} (x_{j})$ is primitive in $R$.

Assume for the rest of the Lemma that  $R = \toba(V)$.

(c). In the situation of (a), resp. (b),
$x_{i}^N = 0$, resp. $\ad x_{i}^r (x_{j}) = 0$.

(d). If $R = \toba (V)$ is an integral domain, then $q_{ii} = 1$
or it is not a root of 1, for all $i$.

(e). If $i\neq j$, then $\ad_{c}(x_{i})^{r} (x_{j}) = 0$
if and only if  $(r)!_{q_{ii}} \prod_{0\le k \le r - 1}
\left(1 - q_{ii}^{k} q_{ij}q_{ji} \right) = 0$.
\end{lema}

\pf (a) and (b) are consequences of the quantum binomial formula,
see {\it e. g.} \cite[Appendix]{AS2} for (b). Then (c) and (d)
follow; the second statement in (c) is also a consequence of (e).
Part (e) is from \cite[Lemma 14]{Ro3}; it can also be shown using
skew-derivations as in \cite[Lemma 3.7]{AS5}. \epf

We now recall a variation of a well-known result of Reshetikhin
on twisting
\cite{Re}.
Let  $({\widehat V}, {\widehat c})$ be another braided vector
space of the same dimension
as $V$, such that
the braiding ${\widehat c}$ is diagonal with matrix
$({\widehat q}_{ij})$ with respect to
a basis ${\widehat x}_{1}, \dots {\widehat x}_{\theta}$.
We define characters ${\widehat \chi}_{1}, \dots, {\widehat \chi}_{\theta}$
of $\Gamma$ by
$$
{\widehat \chi}_{i} (g_{j}) = {\widehat q}_{ji}, \qquad 1\le i,j\le \theta.
$$
We consider ${\widehat V}$ as a Yetter-Drinfeld module
over $\Gamma$ by defining
${\widehat x}_{i} \in {\widehat V}^{{\widehat \chi}_{i}}_{g_{i}}$,
for all $i$.

\begin{prop}\label{reshe}
Assume that for all  $i, j$, $q_{ii} = {\widehat q}_{ii}$  and
\begin{equation}\label{twistequiv}
q_{ij}q_{ji} = {\widehat q}_{ij}{\widehat q}_{ji}.  \end{equation}

Then there exists an $\N$-graded isomorphism of $\ku\Gamma$-comodules
$\psi: \toba(V) \to \toba ({\widehat V})$ such that
\begin{equation}
\psi(x_{i}) = {\widehat x}_{i}, \qquad 1\le i \le \theta.
\end{equation}
Let $\sigma: \Gamma \times \Gamma \to \ku^{\times}$ be the unique bilinear form such that
$\sigma(g_{i}, g_{j}) = {\widehat q}_{ij}  q^{-1}_{ij}$, if $i\le j$,
and is equal to 1 otherwise;
$\sigma$ is a group 2-cocycle and we have, for all $g,h\in \Gamma$,
\begin{align}
\label{psimult}
\psi(xy) &= \sigma(g, h)\psi(x)\psi(y), \qquad x \in \toba (V)_{g},
\quad y \in \toba (V)_{h}; \\ \label{psiad}
 \psi([x, y]_{c}) &=  \sigma(g, h) [\psi(x), \psi(y)]_{c},
\qquad x \in \toba (V)^{\chi}_{g}, \quad y \in \toba (V)_{h}^{\eta}.
\end{align}
\end{prop}
\pf
See \cite[Prop. 3.9 and  Remark 3.10]{AS5}. \epf

\begin{obs}\label{twisting}
In the situation of the proposition, we say that $\toba(V)$ and
$\toba ({\widehat V})$ are \emph{twist-equivalent}; note  that $\toba(V)$
is  twist-equivalent to a $\toba ({\widehat V})$ with
${\widehat q}_{ij} = {\widehat q}_{ji}$ for all $i$ and $j$, since
all the $q_{ij}q_{ji}$'s have  square roots in $\ku$.
\end{obs}

\begin{lema}\label{symmetrizable} Assume that the braiding
with matrix $(q_{ij})$ is generic  and of Cartan type
with generalized Cartan matrix $(a_{ij})$.
Then $(a_{ij})$ is symmetrizable.
\end{lema}

\pf By \cite[Ex. 2.1]{K}, it is enough to show that
$a_{i_{1}i_{2}} a_{i_{2}i_{3}} \dots a_{i_{t-1}i_{t}}a_{i_{t}i_{1}}
= a_{i_{2}i_{1}} a_{i_{3}i_{2}} \dots a_{i_{t}i_{t-1}}a_{i_{1}i_{t}}$
for all $i_1, i_2, \dots, i_t$. But
$$
q_{i_{1}i_{1}}^{a_{i_{1}i_{2}} a_{i_{2}i_{3}}
\dots a_{i_{t-1}i_{t}}a_{i_{t}i_{1}}}
= q_{i_2 i_2}^{a_{i_{2}i_{1}} a_{i_{2}i_{3}}
\dots a_{i_{t-1}i_{t}}a_{i_{t}i_{1}}}
= \dots
= q_{i_{1}i_{1}}^{a_{i_{2}i_{1}} a_{i_{3}i_{2}}
\dots a_{i_{t}i_{t-1}}a_{i_{1}i_{t}}},
$$
by substituting $q_{i_{1}i_{1}}^{a_{i_{1}i_{2}}}
= q_{i_2i_2}^{a_{i_{2}i_{1}}}$, then
$q_{i_2i_2}^{a_{i_{2}i_{1}}} = q_{i_3i_3}^{a_{i_{3}i_{2}}}$ and so on.
The claim follows because $q_{i_{1}i_{1}}$ is not a root of one.
\epf

The following result is due to Rosso, who sketched an argument in \cite[Th. 2.1]{Ro3}.
We include a proof for completeness.

\begin{lema}\label{rossito} \cite{Ro3}.
Let $\ku = \C$. Assume that the braiding
with matrix $(q_{ij})$ is positive  and of Cartan type
with generalized Cartan matrix $(a_{ij})$. Then
$(a_{ij})$ is symmetrizable,  with symmetrizing diagonal matrix
$(d_i)$; and there is a collection of positive  numbers
$(q_I)_{I\in \mathcal X}$ such that $(q_{ij})$
is twist-equivalent to  $({\widehat q}_{ij})$, where
$${\widehat q}_{ij} = q_I^{d_{i}a_{ij}} \text{ for all } i,j \in I.$$
That is, the braiding associated to $({\widehat q}_{ij})$
is of DJ-type.
\end{lema}

\pf We can assume that the braiding is indecomposable; write
$I = \{1, \dots, \theta\}$. By Remark \ref{twisting},
we can assume that $q_{ij} = q_{ji}$, for all $i, j \in I$.
Given $j \in I$,
there exists a sequence
$i_1 =1, i_2, \dots, i_t = j$ of elements in $I$,
such that $a_{i_{\ell}i_{\ell +1}}\neq 0$ for all $\ell$, $1\le \ell < t$.
Then
$$
q_{11}^{a_{i_{1}i_{2}} a_{i_{2}i_{3}} \dots a_{i_{t-1}i_{t}}}
= q_{j j}^{a_{i_{2}i_{1}} a_{i_{3}i_{2}} \dots a_{i_{t}i_{t-1}}},
$$
as in the proof of the previous Lemma. Since
$\alpha_j := a_{i_{1}i_{2}} a_{i_{2}i_{3}} \dots a_{i_{t-1}i_{t}}$
and $\beta_j := a_{i_{2}i_{1}} a_{i_{3}i_{2}} \dots a_{i_{t}i_{t-1}}$
are integers of the same sign, we can find $b\in \N$ and
a family $(d_i)_{i\in I}$ of positive integers such that
$\dfrac{\alpha_j}{\beta_j} = \dfrac{d_j}{b}$. Let $q_I$ be the unique
positive number such that $q_I^b = q_{11}$. Then
$q_{jj}^{\beta_j} = q_{11}^{\alpha_j} = q_{11}^{\beta_j \frac{d_j}{b}}
= q_{I}^{d_j\beta_j}$. Thus $q_{jj}= q_{I}^{d_j}$; and for all $i,j\in I$,
$q_{ij}q_{ji} = q_{I}^{d_i a_{ij}} = q_{I}^{d_j a_{ji}}$,
and the claim follows.
\epf

\begin{obs} The diagonal braiding with matrix
$$
\begin{pmatrix} q & q^{-1} \\ q^{-1} & - q
\end{pmatrix},
$$
where $q$ is not a root of one,
is generic of Cartan type but not of DJ-type.
That is, the preceding Lemma can not be generalized to the
generic case.\end{obs}

\bigbreak
We now state a very elegant description of Nichols
algebras used by Lusztig in a fundamental way \cite{L3}.

\begin{prop}\label{radical-nichols}
Let $(V, c)$ be as above and assume that $q_{ij} = q_{ji}$
for all $i, j$. Let $B_1, \dots, B_{\theta}$ be non-zero elements in
$\ku$. There is a unique  bilinear form
$(\, \vert \, ): T(V) \times T(V) \to \ku$ such that $(1\vert 1) = 1$ and
\begin{flalign}\label{bil1} &(x_j \vert x_j) = \delta_{ij} B_{i}, \qquad \text{for all }
i, j;&\\
\label{bil2} &(x \vert yy') = (x_{(1)} \vert y)(x_{(2)} \vert y'), \qquad \text{for all }
x, y, y'\in T(V);& \\
\label{bil3}
&(xx' \vert y) = (x \vert y_{(1)})(x' \vert y_{(2)}), \qquad \text{for all }
x, x', y\in T(V).&
\end{flalign}
This form is symmetric and also satisfies
\begin{equation}
(x \vert y) = 0, \qquad \text{for all }x\in T(V)_g, y\in T(V)_h, g\neq h \in \Gamma.
\end{equation}
The homogeneous components of $T(V)$ with respect to its usual $\N$-grading are
also orthogonal with respect to $(\, \vert \, )$.

The quotient $T(V) / I(V)$, where  $I(V) = \{x\in T(V): (x \vert
y) = 0 \forall y\in T(V)\}$ is the radical of the form, is
canonically isomorphic to the Nichols algebra of $V$. Thus, $(\,
\vert \, )$ induces a non-degenerate bilinear form on $\toba(V)$,
which will  again be denoted by $(\, \vert \, )$.
\end{prop}

\pf The existence and uniqueness of the form, and the claims about symmetry and
orthogonality, are proved exactly as in \cite[1.2.3]{L3}. It follows from
the properties of the form that $I(V)$ is a Hopf ideal. We now check that
$T(V) / I(V)$ is the Nichols algebra of $V$; it is enough to verify that
the primitive elements of $T(V) / I(V)$ are in $V$. Let $x$ be a primitive element
in $T(V) / I(V)$, homogeneous of degree $n\ge 2$. Then $(x\vert yy') =0$
for all $y$, $y'$ homogeneous of degrees $m, m'\ge 1$ with $m + m' = n$; thus $x = 0$.
\epf

\subsection*{Nichols algebras arising from quantum groups}
Let $(a_{ij})_{1\le i, j\le \theta}$ be a generalized symmetrizable Cartan
matrix \cite{K}; let $(d_{1}, \dots, d_{\theta})$ be positive integers such that
$d_{i}a_{ij} = d_{j}a_{ji}$.
Let $q\in \ku$, $q\neq 0, 1$, and not a root of 1.
We assume that the braided vector space $(V, c)$ is given by the matrix
$q_{ij} = q^{d_{i}a_{ij}}$.

\medbreak
We now want  to derive some precise information about the algebra $\toba(V)$
mainly from \cite{L3}.  We need to consider vector spaces over the field of rational
functions $\kv$.

\bigbreak
Let $(W, d)$ denote a finite-dimensional braided vector
space over $\kv$; we assume that the braiding $d$ is diagonal
with matrix $(v^{d_{i}a_{ij}})$,
with respect to a basis $y_{1}, \dots y_{\theta}$.
We define characters $\eta_{1}, \dots, \eta_{\theta}$ of $\Gamma$ by
$$
\eta_{i} (g_{j}) = v^{d_{i}a_{ij}}, \qquad 1\le i,j\le \theta.
$$
We consider $W$ as a Yetter-Drinfeld module over $\kv \Gamma$ by defining
$y_{i} \in W^{\eta_{i}}_{g_{i}}$, for all $i$.

\bigbreak
Let $v_i := v^{d_{i}}$, $1\le i \le \theta$. We now take
$B_i := \left(1 - v_i^{-2}\right)^{-1} \in \kv$ as in \cite[1.2.3]{L3}.
By proposition \ref{radical-nichols}, $\toba(W)$ is Lusztig's braided
Hopf algebra $\bf f$ (with Cartan datum given by $i . j = d_ia_{ij}$, $\otraavez$.

\bigbreak
Let $\ac := \Q [v, v^{-1}]$; let
$[n]_i := \dfrac{v_i^n - v_i^{-n}}{v_i  - v_i^{-1}}$,
$[r]_i! = [1]_i[2]_i \dots [r]_i$.
Let $\toba(W)_{\ac}$ be the $\ac$-subalgebra of
$\toba(W)$ generated by all
$$
y_i^{(r)} := \frac{y_i^r}{[r]_i!}, \qquad 1\le i \le \theta, \quad r \ge 0.
$$

The canonical bilinear form $(\, \vert \, ): \toba(W) \times \toba(W) \to \kv$
does not restrict to an $\ac$-bilinear form $\toba(W)_{\ac} \times \toba(W)_{\ac}
\to \ac$, since by \cite[1.4.4]{L3},
$$
(y_i^{(r)} \vert y_i^{(r)}) =  v_i^{\frac{r(r + 1)}{2}}
\frac{\left(v_i  - v_i^{-1}\right)^{-r}}{[r]_i!}.
$$
Following an idea of M\"uller \cite{Mu}, we define
$\wtoba(W)_{\ac}$ as the $\ac$-subalgebra of $\toba(W)$ generated by all
$$\widetilde y_i :=  \left(1 - v_i^{-2}\right) y_i = B_i^{-1} y_i, \qquad
1\le i \le \theta.$$

Then $(\, \vert \, )$ restricts to an $\ac$-bilinear form
$(\, \vert \, ): \toba(W)_{\ac} \times \wtoba(W)_{\ac}
\to \ac$,  by the argument in \cite[Lemma 2.2 (a)]{Mu}.

Note that $\toba(W)_{\ac}$ and $\wtoba(W)_{\ac}$ inherit the
Hopf algebra structure from $\toba(W)$, since $\ac$ is a principal ideal domain
with quotient field $\kv$.

\bigbreak
Let $\mathcal W$ be the Weyl group of the generalized
Cartan matrix $(a_{ij})$
\cite[2.1.1]{L3}; then $\mathcal W$ is finite
if and only if $(a_{ij})$ is of finite type.
Let $w \in \mathcal W$ be an element with reduced expression
$w = s_{i_1}s_{i_2} \dots s_{i_P}$, $P \in \N$.
For any $c = (c_1, \dots, c_P) \in \N^P$, let
$$
L(c) := y_{i_1}^{(c_1)} \, T_{i_1}\left(y_{i_2}^{(c_2)} \right)
\dots T_{i_1} \, T_{i_2} \dots T_{i_{P-1}}\left(y_{i_P}^{(c_P)} \right),
$$
where $T_{i_\ell}$ are the $\kv$-algebra automorphisms of $\bf U$ named
$T'_{i_\ell, -1}$ in \cite[37.1.3]{L3}.
Note that $L(c) = L({\bf h}, c, 0, 1)$, with ${\bf h} := (i_1, \dots, i_P)$,
in  \cite[38.2.3]{L3}.  By \cite[41.1.3]{L3},
$$
T_{i_1} \, T_{i_2} \dots T_{i_{\ell-1}}\left(y_{i_\ell}^{(r)} \right)
\in \toba(W)_{\ac}, \qquad 1 \le \ell < P, \quad r\ge 0,
$$
where we have identified $\toba(W)$ with $\bf U^+$ by \cite[3.2.6]{L3},
that is we identify $E_i := \theta_i^+$ in Lusztig's notation with $y_i$.
We denote
$$
z_{\ell} := T_{i_1} \, T_{i_2} \dots T_{i_{\ell-1}}\left(y_{i_\ell}  \right),
\qquad 1 \le \ell \le P.
$$
Then $z_{\ell} \in \toba(W)_{\ac}$ and
$$
z_{\ell}^r = {[r]_{i_{\ell}}!}\,
T_{i_1} \, T_{i_2} \dots T_{i_{\ell-1}}\left(y_{i_\ell}^{(r)} \right),
\qquad 1 \le \ell \le P, \quad r\ge 0.
$$

\begin{teo} (Lusztig). For all $c = (c_1, \dots, c_P), c' \in \N^P$
\begin{equation}\label{lusztig}
(L(c) \vert L(c')) = \delta_{c, c'} \prod_{1\le s \le P}
\prod_{1\le t \le c_s} \left(1 - v_i^{-2t}\right)^{-1}.
\end{equation}
\end{teo}
\pf This follows from \cite[38.2.3 and 1.4.4]{L3}.
\epf

\bigbreak
We regard $\ku$ as an $\ac$-algebra via the algebra map
$\varphi: \ac \to \ku$ given by $\varphi(v) = q$.
We define
$$
\toba(W)_{\ku} := \toba(W)_{\ac}\otimes _{\ac} \ku, \qquad
\wtoba(W)_{\ku} := \wtoba(W)_{\ac}\otimes _{\ac} \ku.
$$
Then $\tobak$, $\wtobak$ are graded braided Hopf algebras
in $\ydg$ and by tensoring with $\ku$ over $\ac$, we get
a bilinear form
$$
(\, \vert \, ): \tobak \times \wtobak \to \ku.
$$
Define $\pi: T(V) \to \tobak$ and $\widetilde\pi: T(V) \to \wtobak$
by $\pi(x_i):= y_i \otimes 1$,
$\widetilde\pi(x_i):= \widetilde y_i \otimes 1$,
$1\le i \le \theta$.
Since $\pi$ and $\widetilde\pi$ induce isomorphisms of braided
vector spaces of $(V, c)$ with $\tobak(1)$ and $\wtobak(1)$,
the composition
$$\begin{CD}
(\, \vert \, ): T(V) \times T(V) @>\pi \times \widetilde\pi>>
\tobak \times \wtobak @>(\, \vert \, )>> \ku
\end{CD}
$$
is the canonical bilinear form of $T(V)$ as in
Proposition \ref{radical-nichols} with scalars $B_i = 1$, $\otravez$.
These arguments allow to adapt many results of Lusztig to the
case when $q$ is not a root of 1.

\medbreak

\begin{teo}\label{gkfinite}\cite[Theorem 15]{Ro3}; \cite[Section 37]{L3}.
Let $(V,c)$ be a braided vector space
of DJ-type.
Then
$\toba (V) \simeq \ku \langle x_1, \dots, x_{\theta} \vert
\ad_c(x_i)^{1-a_{ij}} = 0, \otrvez \rangle$.  \qed \end{teo}

The following Theorem is part of the folklore of quantum groups.

\begin{teo}\label{gkforqgrps} Let $(V,c)$ be a braided vector space
of DJ-type, with generalized Cartan matrix $(a_{ij})$.

(i). If the Gelfand-Kirillov dimension of $\toba (V)$ is finite,
then $(a_{ij})$ is a Cartan matrix of finite type \cite{K}.

(ii). If $(a_{ij})$ is a finite Cartan matrix, then
 the
Gelfand-Kirillov dimension of $\toba (V)$ is finite
and equal to the number of positive roots. \end{teo}

\pf (i). We can assume that the braiding is connected. Let $w \in
\mathcal W$ be an element with reduced expression $w =
s_{i_1}s_{i_2} \dots s_{i_P}$, $P \in \N$. We keep the notation
above. For all $\ell$, $1\le \ell \le P$, we choose an element
$t_{\ell} \in T(V)$ with $\pi (t_{\ell}) = z_{\ell} \otimes 1$,
and set $b_{\ell} :=$  image of $t_{\ell}$ in $\toba (V)$ by the
canonical projection. We claim that the ordered monomials
$b_1^{c_1}b_2^{c_2} \dots b_P^{c_P}$, $c_1, \dots, c_P\ge 0$, are
linearly independent.

\medbreak
To prove the claim, we choose $a_c \in \ac$ such that
$a_cL(c) \in \wtobaka$ for all $c = (c_1, \dots, c_P) \in \N^P$.
If $\wl(c) := a_cL(c) \prod_{1\le s \le P}\prod_{1\le t \le c_s}
 \left(1 - v_i^{-2t}\right) \in \wtobaka$, then
\begin{equation}
(L(c) \vert \wl(c')) = \delta_{c, c'} a_c,
\end{equation}
for all $c, c' \in \N^P$, by \eqref{lusztig}.
It follows that $\pi(t_1^{c_1}t_2^{c_2}   \dots  t_P^{c_P})
= \alpha_c L(c) \otimes 1 \in \wtobak$,
where $\alpha_c = \varphi \left(\prod_{1\le s \le P}
({[c_s]_{i_{s}}!})^{-1}\right)$ is a non-zero scalar in $\ku$.
Choose elements $\wll(c) \in T(V)$ with
$\widetilde\pi(\wll(c)) = \wl(c) \otimes 1$, $c \in \N^P$.
Then for all $c, c' \in \N^P$, we conclude that
$(t_1^{c_1}t_2^{c_2}   \dots  t_P^{c_P} \vert \wll(c'))
= 0$, if $c\neq c'$, and
$(t_1^{c_1}t_2^{c_2}   \dots  t_P^{c_P} \vert \wll(c))
\neq 0$. Since the form $(\quad \vert \quad )$
factorizes over $\toba(V) \times \toba(V)$, the elements
$b_1^{c_1}b_2^{c_2}   \dots  b_P^{c_P}$, $c_1, \dots, c_P\ge 0$,
are linearly independent.

\medbreak
The following statement is well-known:

Let $A$ be a $k$-algebra, where $k$ is a field. Let $a_1, \dots, a_P \in A$
such that the ordered monomials $a_1^{c_1}a_2^{c_2}   \dots  a_P^{c_P}$,
$c_1, \dots, c_P \ge 0$, are linearly independent. Then
the Gelfand-Kirillov dimension of $A$ is $ \ge P$.

We conclude that the Gelfand-Kirillov dimension of $\toba (V)$ is
$ \ge P$.
If the Cartan matrix $(a_{ij})$ is not finite, then
there are elements in the Weyl group of arbitrary length,
and (i) follows.

\medbreak
(ii). Now the Weyl group is finite, and we take the longest element $w_0$.
Let us consider the natural map
$\varphi: \ku \langle x_1, \dots, x_{\theta} \vert
\ad_c(x_i)^{1-a_{ij}}, \otrvez \rangle \to \toba(V)$.
Since the elements
$b_1^{c_1}b_2^{c_2}   \dots  b_P^{c_P}$, $c_1, \dots, c_P\ge 0$,
are the image of the PBW-basis of $\ku \langle x_1, \dots, x_{\theta} \vert
\ad_c(x_i)^{1-a_{ij}}, \otrvez \rangle$, see {\it e. ~g.}
\cite[Proposition 1.7]{deck}, the claim follows.
\epf

\subsection*{Rosso's characterization of Nichols algebras}

We recall some important results of Rosso.

\begin{teo}\label{rosso-qsr} \cite[Lemma 19]{Ro3}
Let $(V,c)$ be a braided vector space of diagonal type.
If $\toba(V)$ has finite Gelfand-Kirillov dimension, then for all $i\neq j$,
there exists $r > 0$ such
that $\ad_{c}(x_{i})^{r} (x_{j}) = 0$.
 \qed
\end{teo}

\medbreak
\begin{cor}\label{rosso-cor}
Let $(V,c)$ be a braided vector space of diagonal type
with indecomposable matrix. Assume that
$\toba(V)$ has finite Gelfand-Kirillov dimension.

(a). If there exists $i$ such that $q_{ii} = 1$, then
$\theta = 1$.

(b). If the braiding is generic then it is of Cartan type.
\end{cor}
\pf This follows  from Theorem \ref{rosso-qsr}
and Lemma \ref{primitivos} (e). \epf

\begin{teo}\label{rosso}  \cite[Theorem 21]{Ro3}
Let $(V, c)$ be a finite-dimensional braided vector space
with positive braiding.
Then   the following are equivalent:

(a). $\toba(V)$ has finite Gelfand-Kirillov dimension.

(b). $(V, c)$ is twist-equivalent to a braiding of DJ-type
with finite Cartan matrix.
\end{teo}

\pf We can assume that the matrix $(q_{ij})$ is indecomposable.
(b) $\implies$ (a) follows from Theorem \ref{gkforqgrps} (ii).
(a) $\implies$ (b).  $(V, c)$ is of Cartan type
by Corollary  \ref{rosso-cor} (b)
with Cartan matrix $(a_{ij})$.
We know that $(a_{ij})$ is symmetrizable by Lemma \ref{symmetrizable}
and that
$(V, c)$ is twist-equivalent to a braiding of DJ-type
by Lemma \ref{rossito}.
By Theorem \ref{gkforqgrps} (i),
the Cartan matrix $(a_{ij})$ is  finite. \epf

\begin{obs} See \cite{AS2} for the analogous problem of characterizing
finite-dimensional braided vector spaces
with diagonal braiding such that  $\toba(V)$ has finite dimension.

\end{obs}

\section{Coradically graded coalgebras}

In this Section we prove a general criterion to determine
the coradical filtration of certain Hopf algebras.
We generalize a method of Takeuchi \cite{T},
who computed the coradical filtration of
$U_q(g)$ in this way; see also \cite{Mu}.
We first extend the definition of coradically graded coalgebras
\cite{CM}.

\medbreak
Let $T\ge 1$ be a natural number. If ${\bf i} = (i_1, \dots, i_T) \in \N^T$,
then we set $\vert {\bf i} \vert = i_1 +  \dots + i_T $.

\begin{defi} An \emph{ $\N^T$-graded coalgebra} is a coalgebra $C$ provided
with an $\N^T$-grading $C = \oplus_{{\bf i} \in \N^T} C({\bf i})$
such that $\Delta C({\bf i})
\subset \oplus_{{\bf j}} C({\bf j}) \otimes C({\bf i - j})$.
An  $\N^T$-graded coalgebra $C$ is  \emph{coradically graded}
if the $n$-th term of the coradical filtration is
$$C_n = \oplus_{{\bf i} \in \N^T, \vert {\bf i} \vert \le n} C({\bf i}),
\quad \forall n\in \N.$$

We denote by $\pi_{\bf i}: C \to C({\bf i})$ the projection associated to the grading.

An  $\N^T$-graded coalgebra $C$ is  \emph{strictly coradically graded}
if $C(0) = C_0$, the coradical of $C$,
and $\Delta_{\bf i,j} : C({\bf i + j} ) \to C({\bf i}) \otimes C({\bf j})$,
$\Delta_{\bf i,j} =  (\pi_{\bf i} \otimes \pi_{\bf j}) \circ \Delta$,
is injective, for all $\bf i,j \in \N^T$.
\end{defi}

\begin{lema}\label{lemacorfilt}
(a).  Let $C$ be a strictly coradically
$\N^T$-graded coalgebra and let $D$ be a strictly coradically
$\N^S$-graded coalgebra.  Then $C \otimes D$
is strictly coradically  $\N^{T + S}$-graded
with respect to the tensor product grading.

(b). If $C$ is strictly
coradically $\N^T$-graded, then it is coradically graded.

(c). If $C$ is coradically $\N$-graded,
then it is strictly coradically graded.\end{lema}

\pf (a)  follows from the definition. We prove (b) by induction on
$n$, the case $n=0$ being part of the hypothesis. Assume $n>0$. If
$c \in \oplus_{\vert {\bf i} \vert \le n} C({\bf i})$, then
$\Delta (c) \in \sum_{\vert {\bf n} \vert = n}C({\bf 0}) \otimes
C({\bf n}) + \oplus_{\vert {\bf j} \vert < n} C({\bf n - j})
\otimes C({\bf j}) \subset C_0 \otimes C + C \otimes C_{n -1}$;
hence $c \in C_n$ (the argument does not need the second
hypothesis). Conversely, let $c \in C_n \cap C({\bf i})$ with
$\vert{\bf i} \vert > n$. By the recursive hypothesis, $\Delta (c)
\in C({\bf 0}) \otimes C + \oplus_{\vert {\bf j} \vert < n} C({\bf
i - j}) \otimes C({\bf j})$. If   ${\bf e} \in N^T$ has $\vert
{\bf e} \vert = 1$, then $\Delta_{\bf e, i-e} (c) = 0$, thus
$c=0$.

(c). See for example \cite[Lemma 2.3]{MiS}.
\epf

We now consider the following situation.
Let $U$ be a Hopf algebra, $H$  a Hopf subalgebra and $N_1, \dots, N_T$
subalgebras of $U$, such that the multiplication induces a linear isomorphism
$$\mu: N_1 \otimes \dots \otimes N_T \otimes H \to U.$$
We assume that the following condition holds:

For all $l$, $1\le l \le T$, $N_lH$ is a Hopf subalgebra of $U$.
Furthermore, the projection
$\pi_l:N_lH \to H$, defined by $\pi_l (nh) = \varepsilon(n)h$, $n \in N_l$, $h \in H$,
is a Hopf algebra map, and $N_l = (N_lH)^{\text{co} \pi_l}$.

Then $N_l$ is a braided Hopf algebra in $\ydhh$ and $N_lH \simeq N_l\# H$.
Let $\Delta_l$ be the comultiplication of the braided Hopf algebra $N_l$
and let $j_l: N_lH  \otimes N_l \to  N_lH  \otimes N_l$ be the map
given by
$j_l(rh \otimes s) = r s_{(-1)} h  \otimes s_{(0)};$
we know that
\begin{equation}\label{twisted-delta}
\Delta (s) = j_l \Delta_l(s), \text{ for all } s\in N_l.
\end{equation}
\begin{lema}\label{lemacorfilt-tecnical} The map
$j:  N_1 \otimes   N_1 \otimes   N_2 \otimes   N_2 \otimes \dots \otimes N_T
\otimes N_T \otimes H \otimes H \to U \otimes U$  given by
\begin{equation}\label{j}
j (u_1 \otimes  v_1 \otimes  u_2 \otimes  v_2  \otimes\dots \otimes
u_T \otimes v_T \otimes h \otimes k) =
j_1(u_1 \otimes  v_1)j_2( u_2 \otimes  v_2) \dots j_T(u_T \otimes v_T)(h \otimes k)
\end{equation}
(product in $U\otimes U$), is a linear isomorphism, and  the following diagram commutes:
\begin{equation}\label{digram-of-j}
\begin{CD}
 N_1 \otimes \dots \otimes N_T \otimes H @> \mu >> U\\
@V \Delta_1 \otimes \Delta_2 \otimes \dots  \otimes\Delta_H VV    @VV\Delta V \\
 N_1 \otimes   N_1 \otimes   N_2 \otimes   N_2 \otimes \dots \otimes N_T
\otimes N_T \otimes H \otimes H @>j>> U \otimes U.
\end{CD}
\end{equation}   \end{lema}

\pf
The commutativity of \eqref{digram-of-j} follows from  \eqref{twisted-delta} since
$\Delta$
is multiplicative. We prove now that $j$ is a linear isomorphism. Let
$z = u_1 \otimes  v_1 \otimes  u_2 \otimes  v_2  \otimes\dots \otimes
u_T \otimes v_T \otimes h \otimes k \in  N_1 \otimes
N_1 \otimes   N_2 \otimes   N_2 \otimes \dots \otimes N_T
\otimes N_T \otimes H \otimes H$. Then
\begin{align*}
j (z) &=
u_1 (v_1)_{(-1)}   u_2  (v_2)_{(-1)} \dots u_T (v_T)_{(-1)}  h
\otimes (v_1)_{(0)}  (v_2)_{(0)} \dots (v_T)_{(0)} k \\
&= u_1 \, (\ad (v_1)_{(-T)}   u_2) \, (\ad ((v_1)_{(-T + 1)} (v_2)_{(-T + 1)})  u_3) \,
\dots \,(\ad ((v_1)_{(-2)} (v_2)_{(-2)}
\dots (v_{T-1})_{(-2)})  u_T) \, \times
\\ &\times  (v_1)_{(-1)}  (v_2)_{(-1)}
\dots (v_T)_{(-1)}  h
\otimes (v_1)_{(0)}  (v_2)_{(0)} \dots (v_T))_{(0)} k
\\
&=  (\mu \otimes  \mu) \left( u_1 \otimes (v_1)_{(-2)} \cdot
\left( u_2 \otimes (v_2)_{(-2)}\cdot  \left( \dots  \otimes (v_{T-1})_{(-2)} \cdot u_T
\right)\right) \otimes  (v_1)_{(-1)}  (v_2)_{(-1)}
\dots (v_T)_{(-1)}  h \right.
\\ & \left. \otimes (v_1)_{(0)}  \otimes   (v_2)_{(0)}
\otimes \dots  \otimes (v_T)_{(0)}
 \otimes k\right).
\end{align*}
Here, in the first equality we applied the rule $ (v_i)_{(-1)}   u_j =
 (\ad (v_i)_{(-2)}   u_j)  (v_i)_{(-1)}$ for $i < j$ systematically;
in the second equality, we use $\cdot$ for the tensor product of copies
of the adjoint representation. We now apply successively several
isomorphisms to
this expression. First, if $M$ is any left $H$-comodule and $N$ is any left $H$-module,
the map
$\phi_{N, M}: N \otimes M \to  N \otimes M$, $\phi_M(n \otimes m) = m_{(-1)} \cdot n
\otimes m_{(0)}$
is an isomorphism with inverse   $\phi_M^{-1} (n \otimes m)
=  \Ss^{-1} (m_{(-1)}) \cdot n\otimes m_{(0)}$.
We apply $\id_{ N_1 \otimes \dots \otimes N_T} \otimes \phi^{-1}_{H, N_1
\otimes \dots \otimes N_T} \otimes \id_H$
to $ (\mu \otimes  \mu)^{-1}j (z)$ and get
$$
 u_1 \otimes (v_1)_{(-1)} \cdot
\left( u_2 \otimes (v_2)_{(-1)}\cdot  \left( \dots  \otimes (v_{T-1})_{(-1)} \cdot u_T
\right)\right) \otimes    h
 \otimes (v_1)_{(0)}  \otimes   (v_2)_{(0)}  \otimes \dots  \otimes (v_T)_{(0)}
 \otimes k;
$$
to this expression, we apply
$\id_{ N_1} \otimes \phi^{-1}_{ N_2 \otimes \dots \otimes N_T\otimes H, N_1} \otimes
\id_{ N_2 \otimes \dots \otimes N_T\otimes H}$
and get
$$
 u_1 \otimes
 u_2 \otimes (v_2)_{(-1)}\cdot  \left( \dots  \otimes (v_{T-1})_{(-1)} \cdot u_T
\right) \otimes    h
 \otimes v_1  \otimes   (v_2)_{(0)}  \otimes \dots  \otimes (v_T)_{(0)}
 \otimes k;
$$
iterating this procedure, we obtain $z$. This shows that $j$ is bijective.
\epf

\begin{teo}\label{teocorfilt} We keep the notations above.
We assume that  $N_l= \oplus_{i \in \N}N_l(i)$ is a
coradically graded Hopf algebra in $\ydhh$,
for all $l$, $1\le l \le T$, and that $H$ is cosemisimple.
Let $U({\bf i}) = \mu\left(  N_1(i_1) \otimes \dots \otimes N_T(i_T) \otimes H\right)$,
for all ${\bf i} = (i_1, \dots, i_T) \in \N^T$.
Then $U = \oplus_{{\bf i} \in \N^T} U({\bf i})$  is a
coradically $\N^T$-graded coalgebra.
\end{teo}

\pf The tensor product coalgebra $ N_1 \otimes \dots \otimes N_T \otimes H$
is a strictly coradically
$\N^T$-graded coalgebra by Lemma \ref{lemacorfilt} (a) and (c)
since each $N_l$ is coradically graded and $H$ is cosemisimple.
Since each $N_l(i)$ is a Yetter-Drinfeld submodule of $N_l$, the map
$j$ in Lemma \ref{twisted-delta} is homogeneous.
Hence, it follows from  Lemma \ref{lemacorfilt-tecnical}
that  $U$ is  strictly coradically $\N^T$-graded coalgebra.
Then $U$ is coradically $\N^T$-graded by  Lemma \ref{lemacorfilt} (b).
\epf

\section{A family of pointed Hopf algebras}\label{const-pointed}

In this Section, we fix
\begin{itemize} \item a free abelian group $\Gamma$ of finite rank
$s$,

\medbreak\item a Cartan matrix $(a_{ij})\in \Z^{\theta\times
\theta}$ of finite type \cite{K}; we denote by $(d_{1}, \dots,
d_{\theta})$ a diagonal matrix of positive integers such that
$d_{i}a_{ij} = d_{j} a_{ji}$, which is minimal with this property;

\medbreak\item a  family $(q_{I})_{I\in \mathcal X}$
of elements in $\ku$ which are not roots of 1;

\medbreak\item elements $g_{1}, \dots, g_{\theta} \in
\Gamma$, characters $\chi_{1}, \dots, \chi_{\theta} \in
\widehat{\Gamma}$ such that $\langle \chi_{i}, g_{i}\rangle = q_{I}^{d_i}$
for all $i$, and
\begin{equation}
\label{cartantype}
\langle \chi_{j}, g_{i}\rangle \langle \chi_{i}, g_{j}\rangle =
q_{I}^{d_{i}a_{ij}},
\qquad \text{for all }  1 \le i, j \le \theta, \quad i\in I.
\end{equation}
\end{itemize}

\begin{defi}\label{link-dat} \cite{AS4}.
We say that two vertices $i$ and $j$
{\em are linkable} (or that $i$ {\em is linkable to} $j$) if
\begin{flalign}\label{link0} &i\not\sim j,&\\
\label{link1} &g_{i}g_{j} \neq 1 \text {  and}& \\
\label{link2}
&\chi_{i}\chi_{j} = \varepsilon.&
\end{flalign}

Here, $\varepsilon$ denotes the trivial representation of $\Gamma$. One can easily see,
{\it cf.} \cite{AS4}, that:
\begin{flalign}
\label{link4} & \text{ If $i$ is linkable to $j$, then }
\chi_{i}(g_{j}) = \chi_{j}(g_{i})^{-1}= \chi_{i}(g_{i}) = \chi_{j}(g_{j})^{-1}.&\\
\label{link5} & \text{ If $i$ and $k$, resp. $j$ and $\ell$,
are linkable, then }
a_{ij} = a_{k\ell}, a_{ji} = a_{\ell k}.&\\
\label{link6} & \text{ A vertex $i$
can not be linkable to two different vertices $j$ and $h$}.&
\end{flalign}

A {\em linking datum}  for $\Gamma$, $(a_{ij})$,
$(q_{I})_{I\in \mathcal X}$,
$g_{1}, \dots, g_{\theta}$ and $\chi_{1}, \dots, \chi_{\theta}$ is
a collection $(\lambda_{ij})_{1 \le i < j \le \theta, i \nsim j}$ of elements
in $\{0, 1\}$
such that $\lambda_{ij}$ is arbitrary if $i$ and $j$
are linkable but 0 otherwise. Given a linking datum,
we say that two vertices $i$ and $j$
{\em are linked} if $\lambda_{ij}\neq 0$.

\begin{obs} A detailed investigation of linking data is carried out
in \cite{Di}; see Theorem 4.6 in \emph{loc. cit.} for a characterization in the generic case.
\end{obs}

The collection $\mathcal D = \mathcal D((a_{ij}),(q_{I}), (g_i), (\chi_i),
(\lambda_{ij}))$, where  $(\lambda_{ij})$ is a linking datum, will be called a {\it
generic datum of
finite Cartan type} for $\Gamma$.
If $\ku = \C$, a generic datum of finite Cartan type
will be called  {\it positive} if $q_{I} > 0$,
for all $I \in \mathcal X$.

\bigbreak
Let $\mathcal D'$ be a generic datum of finite Cartan type over a free abelian group
$\Gamma'$ of finite rank, formed by $(a'_{ij})\in \Z^{\theta'\times \theta'}$,
$(q'_{I})_{I\in \mathcal X'}$,
$g'_{1}, \dots, g'_{\theta}$, $\chi'_{1}, \dots, \chi'_{\theta}$
and a linking datum $(\lambda'_{ij})_{1 \le i < j \le \theta, i \nsim j}$.
The data $\mathcal D$ and $\mathcal D'$ are called {\it isomorphic} if $\theta =
\theta'$, and if there exist a group isomorphism $\varphi : \Gamma \to \Gamma'$, a
permutation $\sigma \in \mathbb S_{\theta}$, and elements $0\neq \alpha_i \in \ku$, for
all $1 \leq i \leq \theta$ such that
\begin{itemize} \item
$\varphi(g_{i}) = g'_{\sigma(i)}$, for all $1\le i \le \theta$;
\item $\chi_{i} = \chi'_{\sigma(i)} \varphi$, for all $1\le i, j \le \theta$;
\item $\lambda_{ij} =
\begin{cases}
 \alpha_i \alpha_j \lambda'_{\sigma(i)\sigma(j)},  &\text{if $\sigma(i) < \sigma(j)$}\\
- \alpha_i \alpha_j  \chi_j(g_i) \lambda'_{\sigma(j)\sigma(i)},  &\text{if $\sigma(i)  >
\sigma(j)$}
\end{cases}$
, for all $1 \leq i < j \leq \theta, i \nsim j.$
\end{itemize}
In this case the triple $(\varphi, \sigma, (\alpha_i))$ will be called an {\it
isomorphism} from $\mathcal D$ to $\mathcal D'$.

Note that then $a_{ij} = a'_{\sigma(i)\sigma(j)}$, for all $1 \leq i,j \leq \theta$. Thus
$\sigma$ is an isomorphism of the corresponding Dynkin diagrams. Indeed, for all $i,j$,
$\chi_{j}(g_{i}) = \chi'_{\sigma(i)}(g'_{\sigma(i)})$, hence $\chi_{j}(g_{i})
\chi_{i}(g_{j}) = \chi_{i}(g_{i})^{a_{ij}} = \chi_{i}(g_{i})^{a'_{\sigma(i) \sigma(j)}}$,
and  $a_{ij} = a'_{\sigma(i)\sigma(j)}$, since $\chi_{i}(g_{i})$ is not a root of 1.

Let Isom($\mathcal D, \mathcal D'$) be the set of all isomorphisms from $\mathcal D$ to
$\mathcal D'$.
\end{defi}

\bigbreak
To state the next theorem,
we  follow the conventions of \cite[Section 4]{AS4}.
We assume that the Cartan matrix is
a matrix of blocks corresponding to the connected components;
that is, for each $I\in \mathcal X$, there exist $c_I, d_I$
such that $I = \{j: c_I\le j \le d_I\}$.

Let $\Phi_{I}$, resp. $\Phi_{I}^{+}$, be the root system,
resp. the subset of positive roots,
corresponding to the Cartan matrix
$(a_{ij})_{i, j\in I}$;
then $\Phi = \bigcup_{I\in \mathcal X}\Phi_{I}$, resp. $\Phi^+ =
\bigcup_{I\in \mathcal X}\Phi_{I}^{+}$ is the root system, resp. the subset of positive
roots,
corresponding to the Cartan matrix
$(a_{ij})_{1\le i, j\le \theta}.$
 Let $\alpha_{1}, \dots,\alpha_{\theta}$  be the set of simple roots.

\smallbreak
Let $\mathcal W_I$ be the Weyl group  corresponding to the Cartan matrix
$(a_{ij})_{i, j\in I}$; we identify it with a subgroup of the Weyl group
$\mathcal W$ corresponding to the Cartan matrix
$(a_{ij})$. We fix
a reduced decomposition of the longest element $\omega_{0, I}$ of $\mathcal W_{I}$
in terms of simple reflections.
Then we obtain a reduced decomposition of the longest element $\omega_{0}
= s_{i_1} \dots s_{i_P}$ of $\mathcal W$
from the expression of $\omega_{0}$ as product of the $\omega_{0, I}$'s in some
fixed order of the components, say the order arising from the order of the vertices.
Therefore
$\beta_{j} := s_{i_1} \dots s_{i_{j-1}}(\alpha_{i_j})$
is a numeration of $\Phi^+$.

\bigbreak
We fix a finite-dimensional Yetter-Drinfeld module $V$ over $\Gamma$
with a basis $x_1, \dots, x_{\theta}$ with
$x_{i}\in V^{\chi_i}_{g_i}$, $1\le i \le \theta$.
Note that
\begin{equation}\label{noniso}
V_{g_{i}}^{\chi_{i}} \ncong V_{g_{j}}^{\chi_{j}} \text{ in } \ydg, \text{ for all } 1
\leq i,j \leq \theta, i \neq j;
\end{equation}
see the proof of Lemma \ref{det-basis}.

\bigbreak
Lusztig defined root vectors $X_{\alpha}$,
$\alpha \in \Phi^{+}$ \cite{L3},
in the case of braidings of DJ-type; these are the elements
$b_1, \dots, b_P$ in the proof of Theorem \ref{gkforqgrps};
they can be expressed as iterated braided commutators.
As in \cite{AS4}, this definition can be extended
to generic braidings of finite Cartan type, first
in the tensor algebra $T(V)$, and then in suitable quotients.

We fix a $\Z$-basis $Y_{h}, 1 \leq h \leq s$ of $\Gamma$.

\begin{teo}\label{construction} Let $\mathcal D = \mathcal D((a_{ij}),(q_{I}), (g_i),
(\chi_i),  (\lambda_{ij}))$ be a generic datum of finite Cartan type for the free abelian
group $\Gamma$ of finite rank.
Let $U({\mathcal D})$  be the algebra presented by generators $a_{1}, \dots, a_{\theta}$,
$y_{1}^{\pm 1}, \dots, y_{s}^{\pm 1}$ and relations
\begin{align}\label{relations}
 y_{m}^{\pm 1}y_{h}^{\pm 1} &= y_{h}^{\pm 1}y_{m}^{\pm 1}, \quad y_{m}^{\pm 1}
y_{m}^{\mp 1} = 1, \qquad 1 \le m,h  \le s,\\\label{relations1}
y_{h}a_{j} &= \chi_{j}(Y_{h})a_{j}y_{h}, \qquad 1 \le h \le s,\, 1 \le j \le \theta,\\
\label{relations2}
(\ad a_{i})^{1 - a_{ij}}a_{j} &= 0, \qquad 1 \le i \neq j \le \theta, \quad i \sim j, \\
\label{relations3}
a_{i}a_{j} - \chi_{j}(g_{i})a_{j}a_{i} &= \lambda_{ij}(1 - g_{i}g_{j}),
\qquad 1 \le i < j \le \theta, \quad i \not\sim j;
\end{align}
then $U({\mathcal D})$ is a pointed Hopf algebra
with structure determined  by
\begin{equation}\label{sk-gl}
\Delta y_{h} = y_{h}\otimes y_{h}, \qquad
\Delta a_{i} = a_{i}\otimes 1 + g_{i} \otimes a_{i},
\qquad 1 \le h \le s,\, 1 \le i \le \theta.
\end{equation}

Furthermore,  $U({\mathcal D})$ has a PBW-basis
given by monomials in the root vectors $b_1 :=a_{\beta_1}, \dots,
b_P :=a_{\beta_P}$.
The coradical filtration of $U({\mathcal D})$ is given by
\begin{equation}
U({\mathcal D})_{N} = \text{ span of } a_{i_1}a_{i_2}   \dots
a_{i_r} y: \quad r \le N,  \quad y\in  G(H).
\end{equation}
There is an isomorphism of graded Hopf algebras
$\psi: \toba(V)\# \ku \Gamma \to \gr U({\mathcal D})$,
given by $ x_i \# 1 \mapsto \overline{a_i}$, $1 \le i \le \theta$,
$1 \# Y_{h} \mapsto  y_{h}$, $1 \le h  \le s$.

$U({\mathcal D})$  has finite Gelfand-Kirillov dimension and
is a domain. \end{teo}

In  \eqref{relations3} and \eqref{sk-gl}, the elements $g_i$, $ 1 \leq i \leq \theta$,
must be read as the word in the generators
$y_{h}, 1 \leq h \leq s$  given by
$y_{1}^{t_{i,1}} \cdots y_{s}^{t_{i,s}}$, if
 $g_i = Y_{1}^{t_{i,1}} \cdots Y_{s}^{t_{i,s}}$,
where $t_{i,1}, \cdots, t_{i,s}$ are integers.

The relations \eqref{relations2} are the quantum Serre relations. By \eqref{relations1}
and \eqref{sk-gl}, $(\ad a_i)(a_j) = a_i a_j - \chi_j(g_i) a_j a_i = (\ad_{c}
a_{i})(a_{j})$. Here, $\ad_{c}$ is the braided adjoint representation in the tensor
algebra of $V$. Hence the left hand side of \eqref{relations2} should be more formally
written as
$$ (\ad_{c} a_{i})^{1-a_{ij}}(a_{j}) = \sum_{l=0}^{1-a_{ij}} (-1)^{l}
\binom{1-a_{ij}}{l}_{q_{ii}} q_{ii}^{\frac{l(l-1)}{2}} q_{ij}^{l} a_{i}^{1-a_{ij}-l}
a_{j} a_{i}^{l},$$
where $q_{ij}:= \chi_{j}(g_{i})$ for all $i,j$.

\pf {\it Step I}.
It is not difficult to see that $\Delta$
is well-defined by  \eqref{sk-gl}, using
Lemma \ref{primitivos} (b) for the quantum Serre relations.

\medbreak
 {\it Step II}. As in \cite[Th. 4.5]{AS4},
we deduce from Theorem  \ref{gkfinite}, via
Proposition  \ref{reshe}, that
$\toba (V) \simeq \ku \langle x_1, \dots, x_{\theta} \vert
\ad_c(x_i)^{1-a_{ij}} = 0, \otrvez \rangle$.

\medbreak
{\it Step III}. We now prove the statement
about the PBW-basis. We argue exactly as in
the proof of \cite[Th. 5.17]{AS4}; we proceed
by induction on the number of connected components.
If the Dynkin diagram is connected, the claim follows
from Step II as in \cite[Lemma 5.18]{AS4}.
For the inductive step, one repeats the argument
in  \cite[Lemma 5.19]{AS4}.
We assume there exists $\ttheta < \theta$ such that
$J = \{1, \dots, \ttheta\} \in \mathcal X$.
Let $\Upsilon := <Z_1>\oplus \dots \oplus <Z_{\ttheta}>$,
be a free abelian group of rank $\ttheta$.
Let $\eta_{j}$ be the unique character of $\Upsilon$ such that
$\eta_{j}(Z_{i}) = \chi_{j}(g_{i})$, $1\le i,j \le \ttheta$.

Let $\mathcal{D}_1 $ be the generic datum over
$\Gamma$ given by
$ (a_{ij})_{\ttheta <  i, j\le \theta}$,
$(g_{i})_{\ttheta < i \le \theta}$,
$(\chi_{j})_{\ttheta < j \le \theta}$,
$(\lambda_{ij})_{\ttheta < i < j \le \theta,  i\not\sim j}$.
Let $\bc := U(\mathcal{D}_1)$,
with generators $b_{\ttheta +1}$, \dots, $b_{\theta}$
(instead of the $a_{i}$'s) and  $y_{1}, \dots, y_{s}$.

Let $\mathcal{D}_2$ be the generic datum over
$\Upsilon$ given by $(a_{ij})_{1\le i, j\le \ttheta}$,
$(Z_{i})_{1 \le i \le \ttheta}$,
$(\eta_{j})_{1 \le j \le \ttheta}$,
with empty linking datum.
Let $\uc := U(\mathcal{D}_2)$ with generators $u_{1},
\dots, u_{\ttheta}$ (instead of the $a_{i}$'s)
and  $z_1, \dots, z_{\ttheta}$.

Then the analogue of \cite[Lemma 5.19]{AS4} holds
replacing the dual of $\bc$ by the Hopf dual.
Notice that the argument in {\it loc. cit.} works
since algebra maps and skew-derivations are in the
Hopf dual. The proof is actually easier
because there are less relations to check.
Finally, we form the Hopf algebra
$(\uc \otimes \bc)_{\sigma}$ as in {\it loc. cit.},
and consider the quotient $\acw$ of
$(\uc \otimes \bc)_{\sigma}$  by the central Hopf subalgebra
$\ku[(z_{i} \otimes g_{i}^{-1}): 1\le i \le \ttheta]$.
The same argument as in {\it loc. cit.} shows that $\acw \simeq U({\mathcal D})$
as Hopf algebras. On the other hand,
the monomials
$b_1^{c_1}b_2^{c_2}   \dots  b_P^{c_P} y$,
$c_j \in \N$,  $1\le j \le P$, $y\in \Upsilon\times\Gamma$
form a basis of $(\uc \otimes \bc)_{\sigma}$.
Since $\Upsilon\times\Gamma $ splits as the product
of $\Gamma$ and the group generated by
$(z_{i} \otimes g_{i}^{-1})$, $1\le i \le \ttheta$,
we conclude that the images of the monomials
$b_1^{c_1}b_2^{c_2}   \dots  b_P^{c_P} y$,
$c_j \in \N$,  $1\le j \le P$, $y\in \Gamma$
form a basis of $\acw$.

\medbreak
{\it Step IV}. The claim about the GK-dimension
follows from the previous step. The claim about the
coradical filtration also
follows from the previous step,
together with Theorem \ref{teocorfilt}.
Indeed, Theorem \ref{teocorfilt} applies
since Nichols algebras are coradically graded
by definition.

\medbreak
{\it Step V}. The existence of $\psi$ follows from the
definition of Nichols algebras, since by
the statement about the coradical filtration
we have a monomorphism of Yetter-Drinfeld modules
$V \to  U({\mathcal D})_{1} / U({\mathcal D})_{0}$.
Since the restriction of  $\psi$ to the first
term of the coradical filtration $(\toba(V)\# \ku \Gamma)_1$
is injective, $\psi$ is injective \cite[Th. 5.3.1]{M}.
It is surjective by the PBW-basis claim.
Therefore $\psi$ is bijective.

\medbreak
{\it Step VI}. We finally  prove that $ U({\mathcal D})$
is a domain. Here we follow \cite[Corollary 1.8]{deck}.
We introduce an $\N^{P + 1}$-filtration on $ U({\mathcal D})$ by
the degree defined by
$$\deg \, (b_1^{c_1}b_2^{c_2}   \dots  b_P^{c_P} y) =
(c_1, c_2, \dots, c_P, \sum_{1\le j \le P} c_j \height \beta_j),
\quad, c_j \in \N, \quad 1\le j \le P, \quad y\in \Gamma;$$
here $\height \beta$ is the height of the root $\beta$
as in \cite{deck}. We claim that this is an algebra filtration.
If $\alpha = e_1\alpha_1 +
\dots + e_{\theta}\alpha_{\theta}$, we set
$g_{\alpha} = g_{1}^{e_1} \dots g_{\theta}^{e_{\theta}}$,
$\chi_{\alpha} =
\chi_{1}^{e_1} \dots \chi_{\theta}^{e_{\theta}}$.
Recall that $b_k = a_{\beta_k}$, $1\le k \le P$.
To prove the claim one has to verify that
for all $k > l$
\begin{equation}\label{levsoib}
b_k b_l - \chi_{\beta_l}(g_{\beta_k}) b_l b_k =
\sum_{{\bf c} \in \N^P} \rho_{{\bf c}} \,
b_1^{c_1}b_2^{c_2}   \dots  b_P^{c_P},
\end{equation}
where $\rho_{{\bf c}} \in \ku$, and  $\rho_{{\bf c}} = 0$
unless $\deg \, (b_1^{c_1}b_2^{c_2}   \dots  b_P^{c_P})
< \deg \, (b_l b_k)$.

If $\beta_l$ and $\beta_k$ have support in the
same connected component, then \eqref{levsoib}
follows from the formula of Levendorskii and Soibelman
 \cite[Lemma 1.7]{deck} using Lemma \ref{reshe}.
Note that here we are using the special order of the
set of positive roots, in which the roots
with support on the first component are smaller
than the roots in the second, and so on.

If $\beta_l$ and $\beta_k$ have support in
different connected components, then \eqref{levsoib}
follows from the linking relations since the
monomials in any root vector are homogeneous.

As in \cite[Corollary 1.8]{deck}, it is not difficult to verify
that the associated graded ring is a domain, which implies
that $ U({\mathcal D})$ is a domain.
\epf

To describe the isomorphisms between Hopf algebras $U({\mathcal D})$ and  $U({\mathcal
D'})$, we first formulate a Lemma which is needed in this general form in the proof of
the
main theorem \ref{fingrowth-lifting}.

\begin{lema}\label{A1}
Let $\Gamma$ be an abelian group, $A$ a pointed Hopf algebra with
coradical $A_0 = \ku \Gamma$, $V \in \ydg$ with $\ku$-basis $x_i \in
V_{g_{i}}^{\chi_{i}}, g_{i} \in \Gamma, \chi_{i} \in \VGamma, 1 \leq i \leq \theta.$

Assume that
$\text{gr}A \cong \toba(V) \# \ku \Gamma$ as graded Hopf algebras,
 and
\begin{equation}\label{not1}
\chi_{i} \neq \varepsilon, \text{ for all }1 \leq i \leq \theta.
\end{equation}
Then the first term of the coradical filtration of $A$ is
$$A_1 = A_0 \oplus \bigoplus_{g \in \Gamma, 1 \leq i \leq \theta} {\mathcal
P}_{gg_{i},g}(A)^{\chi_{i}}.$$
If $(g_i, \chi_i) \neq (g_j, \chi_j)$ for all $i \neq j$, then the vector spaces
${\mathcal P}_{gg_{i},g}(A)^{\chi_{i}}$ are one-dimensional for all
$g \in \Gamma, 1 \leq i \leq \theta$.
\end{lema}
\pf
By assumption, $A_1/A_0 \cong V \# \ku \Gamma$, and  the elements
$x_{i} \# g \in{\mathcal P}_{gg_{i},g}(\toba(V) \# \ku \Gamma )^{\chi_{i}}$,
$1\leq i \leq \theta$, $g \in
\Gamma$, $x \in {\mathcal P}_{gg_{i},g}(A)^{\chi_{i}}$
form a $\ku$-basis of $V \# \ku
\Gamma$. Hence
\begin{equation}\label{not2}
A_1/A_0 \cong \bigoplus_{1 \leq i \leq \theta} (A_1/A_0)^{\chi_{i}}.
\end{equation}

We first show that $A_1$ is locally finite under the adjoint
action of $\Gamma$. By the theorem of Taft and Wilson
\cite[Theorem 5.4.1]{M}, $A_1 = A_0 + (\bigoplus_{g,h \in \Gamma}
{\mathcal P}_{g,h}(A))$. Hence it is enough to prove that for $g,
h \in \Gamma$, ${\mathcal P}_{g,h}(A)$  is locally finite.

Since ${\mathcal P}_{g,h}(A) / \ku (g-h)$ is embedded into
$A_1/A_0$, it follows from \eqref{not2} that ${\mathcal
P}_{g,h}(A) / \ku (g-h)$ and  ${\mathcal P}_{g,h}(A)$ are locally
finite.

We claim that $A_1$ is completely reducible as $\Gamma$-module.
Indeed, let $U$ be any locally finite $\Gamma$-module, let $\chi \in \VGamma$, and
let $U^{(\chi)} = \{u\in U: \exists s > 0$ such that $(g - \chi(g))^s (u) = 0 \forall
g\in \Gamma\}$. Then $U = \oplus_{\chi \in \VGamma} U^{(\chi)}$ (see for instance
\cite[Th. 1.3.19]{dixmier}). Now, $A_0 \subset (A_1)^{\varepsilon} \subset
(A_1)^{(\varepsilon)}$,
but by \eqref{not1}, they are all three equal. The claim follows.

Hence
$A_1 = A_0 \oplus \bigoplus_{1 \leq i \leq \theta} (A_1)^{\chi_{i}}$.
By the theorem of Taft and Wilson again
$$A_1 = A_0 \oplus \bigoplus_{g,h \in \Gamma, \varepsilon \neq \chi \in \VGamma}
{\mathcal P}_{g,h}(A)^{\chi},$$
and the lemma follows.
\epf
Note that \eqref{not1} holds for generic braidings, since in this case no $\chi_i(g_i)$
is a root of 1.

If $A,B$ are Hopf algebras, we denote the set of all Hopf algebra isomorphisms from $A$
to $B$ by Isom($A,B$).

\begin{teo}\label{clasif} Let
${\mathcal D}$ and ${\mathcal D'}$ be generic
data of finite Cartan type for $\Gamma$.
Then the Hopf algebras $U({\mathcal D})$
and $U({\mathcal D'})$ are isomorphic
if and only if ${\mathcal D}$ is isomorphic to ${\mathcal D'}$.

More precisely, let $a_1, \cdots, a_{\theta}$ resp. $a'_1, \cdots, a'_{\theta}$
be the simple root vectors in $U({\mathcal D})$ resp. $U({\mathcal D'})$ of Theorem
\ref{construction}, and let $g_1, \cdots, g_{\theta}$ resp. $g'_1, \cdots, g'_{\theta}$
be the group-like elements in ${\mathcal D}$ resp. ${\mathcal D'}$ . Then the map
$$\text{Isom}(U({\mathcal D}),U({\mathcal D'})) \to \text{Isom}({\mathcal D},{\mathcal
D'}),$$
given by $\phi \mapsto (\varphi, \sigma, (\alpha_i))$, where $\varphi(g) = \phi (g),
\varphi(g_i) = g'_{\sigma(i)}, \phi(a_i) = \alpha_i a'_{\sigma(i)}$, for all $g \in
\Gamma, 1\leq i \leq \theta$, is bijective.
\end{teo}
\pf
Let $V$ resp. $V'$  be the Yetter-Drinfeld module of
the infinitesimal braiding of $A:= U({\mathcal D})$ resp. of $A':= U({\mathcal D})$. Let
$\phi : A \to A'$ be an isomorphism of Hopf algebras. Then $\phi$ induces isomorphisms
$A_0 \to A'_0$ and $A_1 \to A'_1$. Hence $\phi$ defines an isomorphism of groups $\varphi
: \Gamma \to \Gamma$, and for all $g,h \in \Gamma, \chi \in \VGamma$, a linear
isomorphism
$$ {\mathcal P}_{g,h}(A)^{\chi} \cong {\mathcal P}_{\varphi(g),\varphi(h)}(A')^{\chi
\varphi^{-1}}.$$
By Theorem \ref{construction}, the assumptions of Lemma \ref{A1} are satisfied for $A,V$
and $A',V'$. Then it follows from Lemma \ref{A1} and \eqref{noniso} that there is a
uniquely determined permutation $\sigma \in \mathbb{S}_{\theta}$ such that $\phi$ induces
an isomorphism
$${\mathcal P}_{g_{i},1}(A)^{\chi_{i}} \cong {\mathcal
P}_{g'_{\sigma(i)},1}(A')^{\chi'_{\sigma(i)}}, \text{ with } \varphi(g_{i}) =
g'_{\sigma(i)}, \chi_{i} \varphi^{-1} = \chi'_{\sigma(i)}, \text{ for all } 1 \leq i \leq
\theta.$$
Moreover, since for all $i$, ${\mathcal P}_{g_{i},1}(A)^{\chi_{i}}$ and ${\mathcal
P}_{g'_{\sigma(i)},1}(A)^{\chi'_{\sigma(i)}}$ are one-dimensional with basis $a_i$ and
$a'_i$, there are non-zero scalars $\alpha_{i} \in \ku$ with $\phi(a_i) = \alpha_i
a'_{\sigma(i)}$, for all $1 \leq i \leq \theta$.

Then the elements $\phi(a_i), 1 \leq i \leq \theta$, satisfy the Serre relations
\eqref{relations2}, and they satisfy \eqref{relations3} if and only if the triple
$(\varphi, \sigma, (\alpha_i))$ is an isomorphism of generic data.
Thus the map $\text{Isom}(U({\mathcal D}),U({\mathcal D'})) \to \text{Isom}({\mathcal
D},{\mathcal D'})$ in the theorem is well-defined and injective. Surjectivity of this map
follows from the description of the Hopf algebras $U({\mathcal D})$ and $U({\mathcal
D'})$ in Theorem \ref{construction}.
\epf

The main reason why the proof of the preceding theorem works
is the knowledge of the
coradical filtration. The same ideas allow to determine all
Hopf subalgebras of $U({\mathcal D})$.

\section{Pointed Hopf algebras with generic braidings}\label{pos-br}
We are going to show that the class of Hopf algebras
described in the previous section has an intrinsic description.

The following key Lemma implies that  pointed Hopf algebras belonging to a natural class
are generated by group-like and skew primitive elements.

\begin{lema}\label{fingrowth-nichols} (a).  Let $S = \oplus_{n\in \N} S(n)$
be a graded braided Hopf algebra such that
$S(0) = \C.1$, $V :=S(1)$ is finite-dimensional
and generates $S$ as an algebra.
Assume that $S$ has finite Gelfand-Kirillov dimension and
that $V$ has
positive  braiding. Then $S$ is a Nichols algebra.

(b). Let $R$ be as $S$ in (a), except that we assume $P(R) = R(1)$ instead of generation
in degree 1.
Then $R$ is a Nichols algebra.
\end{lema}

\pf (a).
 $\toba(V)$  has finite Gelfand-Kirillov dimension
since it is a quotient of $S$. Assume first that the matrix is
indecomposable. We can then apply Theorem \ref{rosso}; let
$(a_{ij})$, $(d_{1}, \dots, d_{\theta})$ and $q$ be such that
$q_{ij}q_{ji} = q^{d_{i}a_{ij}}$ for all $i \neq j$ and $q_{ii} =
q^{d_{ii}}$.

Let $i\neq j$. We claim that $z_{2} = \ad_{c}(x_{i})^{1 - a_{ij}} (x_{j}) = 0$ in $S$.
Indeed, let $z_{1} = x_{i}$, suppose that $z_{2} \neq 0$ and consider the
two-dimensional subspace $W$ of $S$ generated by the primitive
elements $z_{1}$ and $z_{2}$.

We claim  that $\toba (W)$ has finite Gelfand-Kirillov dimension. For,
let $T$ be the subalgebra of $S$ generated by $W$; then
the graded Hopf algebra $\gr (T \# \C \Gamma)$ has finite
Gelfand-Kirillov dimension, and contains $\toba (W)$.

Also, the braiding of $W$ is given by the matrix:

$$
\begin{pmatrix} q_{ii} & q_{ii}^{1-a_{ij}} q_{ij} \\
q_{ii}^{1-a_{ij}} q_{ji} & q_{ii}^{(1-a_{ij})^2} q_{ij}^{1-a_{ij}}q_{ji}^{1-a_{ij}}
q_{jj}
\end{pmatrix} = \begin{pmatrix} q^{d_{i}} & q^{d_{i}(1-a_{ij})} q_{ij} \\
q^{d_{i}(1-a_{ij})} q_{ji} & q^{d_{i} -d_{i}a_{ij} + d_{j}}
\end{pmatrix}
$$
By Theorem \ref{rosso-qsr} and Lemma \ref{primitivos}, there exists $k \ge  0$
such that $1 = q^{d_{i}k + 2d_{i}(1-a_{ij})} q_{ij} q_{ji}$, hence
$0 = d_{i}k + 2d_{i}(1-a_{ij}) + d_{i} a_{ij} = d_{i}(k + 2 - a_{ij})$,
a contradiction. This shows that $z_{2} = 0$.

Therefore, we have an epimorphism of braided graded Hopf algebras
$\toba(V) \to S$, by Step II
of Theorem \ref{construction}, which is the identity in degree 1.
Hence $\toba(V) \simeq S$.

\medbreak Assume now that the matrix is decomposable. Let $i$, $j$
belong to different components; in particular $q_{ij}q_{ji} = 1$.
We claim that $x_{i}x_{j} = q_{ij} x_{j}x_{i}$. If not, let $z_{1}
:= x_{i}$ and $z_{2} :=x_{i}x_{j} - q_{ij} x_{j}x_{i}$, that is
primitive by Lemma \ref{primitivos} (b). Consider as before the
subspace $W$ of $S$ generated by $z_{1}$ and $z_{2}$. As before,
$\toba (W)$ has finite Gelfand-Kirillov dimension. Now the
braiding of $W$ is given by the matrix:

$$ \begin{pmatrix} \widetilde q_{11} &\widetilde  q_{12} \\
\widetilde q_{21} &\widetilde  q_{22} \end{pmatrix} := \begin{pmatrix} q_{ii} & q_{ii}
q_{ij} \\
q_{ii} q_{ji} & q_{ii}   q_{jj} \end{pmatrix}. $$
By Theorem \ref{rosso-qsr} and Lemma \ref{primitivos} again, there exists $k \ge  0$
such that $1 = q_{ii}^{k + 2} q_{ij} q_{ji} = q_{ii}^{k + 2}$, a contradiction.
This concludes the proof of (a).

Finally, (a) and (b) are equivalent by \cite[Lemma 5.5]{AS2}
and the definition of finite Gelfand-Kirillov dimension.
Indeed, we can assume that the homogeneous components
of $R$ are finite-dimensional, by replacing if necessary $R$
by the subalgebra generated by $R(1)$ and any finite-dimensional coalgebra.
Note that braiding of the dual of $V$ is again positive.
\epf

\begin{teo}\label{fingrowth-lifting}
Let $A$ be a pointed Hopf algebra with finitely generated abelian group
$G(A)$, and positive infinitesimal braiding.
Then the following are equivalent:

(a). $A$ is a domain with finite Gelfand-Kirillov dimension.

(b). The group $\Gamma := G(A)$ is free abelian of finite rank,
and there exists a positive datum $\mathcal D$ for $\Gamma$ such
that $A \simeq U({\mathcal D})$ as Hopf algebras.
 \end{teo}

\pf (b) $\implies$ (a): this is Theorem \ref{construction}.

(a) $\implies$ (b). Consider the diagram $R$ of $A$. By
\cite[6.5]{KL}, $\gr A$ has finite GK-dimension; hence both $R$
and $\ku \Gamma$ also have finite GK-dimension. It is clear then
that $\Gamma$ should be a free abelian group of finite rank, say
$s$.  From Theorem \ref{rosso}, Lemma \ref{det-basis} and  Lemma
\ref{fingrowth-nichols} (b), we conclude the existence of the
finite Cartan matrix $(a_{ij})$, the family $(q_{I})_{I\in
\mathcal X}$, and the  elements $g_{1}, \dots, g_{\theta} \in
\Gamma$, $\chi_{1}, \dots, \chi_{\theta} \in \widehat{\Gamma}$
satisfying \eqref{cartantype}, and such that no $q_I$ is a root of
1 (in fact $q_I > 0$ and not 1 for all $I$), $R = \toba(V)$ where
$V \in \ydg$ has a basis $x_i \in V_{g_{i}}^{\chi_{i}}, 1 \leq i
\leq \theta$.

Since $\text{gr}A \cong \toba(V) \# \ku \Gamma$ as graded Hopf algebras, and
$\chi_i(g_i) \neq 1$ for all $1 \leq i \leq \theta$, it follows from Lemma \ref{A1} that
the first term of the coradical filtration of $A$ is
\begin{equation}\label{A11}
A_1 = A_0 \oplus \bigoplus_{g \in \Gamma, 1 \leq i \leq \theta} {\mathcal
P}_{gg_{i},g}(A)^{\chi_{i}}.
\end{equation}

\medbreak We can then choose $a_{i} \in {\mathcal P}_{g_{i},
1}(A)^{\chi_{i}}$ such that the class of $a_{i}$ in $\gr A(1)$
coincides with $x_{i}\# 1$. Let $y_{1}, \cdots, y_{s}$ be free
generators of $G(A)$. It is clear that relations \eqref{relations}
and \eqref{relations1} hold.

Let $i\neq j$. We claim:

\medbreak
(i). There exists no $\ell$, $1\le  \ell
 \le \theta$, such that $g_{i}^{1-a_{ij}}g_{j} = g_{\ell}$,
 $\chi_{i}^{1-a_{ij}}\chi_{j} = \chi_{\ell}$.

\medbreak
(ii). If $i\sim j$, then $\chi_{i}^{1-a_{ij}}\chi_{j} \neq \varepsilon$.

\medbreak
We prove (i). Assume that $g_{i}^{1-a_{ij}}g_{j} = g_{\ell}$,
$\chi_{i}^{1-a_{ij}}\chi_{j} = \chi_{\ell}$ for some $\ell$.
Then
$$q_{I}^{d_{i}a_{i\ell}} = \langle \chi_{\ell}, g_{i}\rangle \langle \chi_{i},
g_{\ell}\rangle =
q_{I}^{2d_{i} (1-a_{ij})}  \langle \chi_{j}, g_{i}\rangle \langle \chi_{i}, g_{j}\rangle=
q_{I}^{d_{i}(2-a_{ij})};$$
we conclude that $2 = a_{ij} + a_{i\ell}$. The only possibility is $a_{ij} = 0$ and
$l=i$. Then $g_j = 1$ which is impossible.

We prove  (ii).
Assume that $\chi_{i}^{1-a_{ij}}\chi_{j} = \varepsilon$, $i\neq j$, $i\sim j$.
Evaluating at $g_i$, we get
$1 = q_{ii}^{1 - a_{ij}} q_{ij} = q_{ii} (q_{ij}q_{ji})^{-1} q_{ij} = q_{ii}
q_{ji}^{-1}$, so that
$q_{ii} = q_{ji}$. Evaluating at $g_j$, we get
$1 = q_{ji}^{1 - a_{ij}} q_{jj}$; hence $q_{jj}= q_{ii} ^{a_{ij}-1}$. Hence
$0  = d_{i}(1 - a_{ij}) + d_{j}$; this is a contradiction.
The claim is proved.

\medbreak We next show that the $a_i$'s satisfy the quantum Serre
relations \eqref{relations2} and the linking relations
\eqref{relations3}. If  $i\neq j$, then $(\ad a_{i})^{1 -
a_{ij}}a_{j} \in {\mathcal
P}_{g_{i}^{1-a_{ij}}g_{j},1}(A)^{\chi_{i}^{1-a_{ij}}\chi_{j}}$ by
Lemma \ref{primitivos} (b). If $i\sim j$, taking into account
\eqref{A11}, (i) and (ii), we see that the quantum Serre relations
\eqref{relations2} hold in $A$.

\medbreak
Finally, assume that $i\not\sim j$; if $0 \neq (\ad a_{i}) a_{j}
\in {\mathcal P}_{g_{i}g_{j},1}(A)^{\chi_{i}\chi_{j}}$,
then $\chi_{i}\chi_{j} = \varepsilon$ by \eqref{A11} and (i).
So that
$a_{i}a_{j} - \chi_{j}(g_{i})a_{j}a_{i}
= \lambda_{ij}(1 - g_{i}g_{j})$
for some $\lambda_{ij}\in \ku$,
where $\lambda_{ij} = 0 $ when $\chi_{i}\chi_{j} \neq \varepsilon$.
But we can also choose $\lambda_{ij} = 0 $ when $g_{i}g_{j} = 1$.
By \eqref{link6}, we can rescale a generator $a_{i}$
with $\lambda_{ij}\neq 0$ to have  $\lambda_{ij} = 1$.
Hence, $(\lambda_{ij})$ is a linking datum for
$(a_{ij})$, $g_{1}, \dots, g_{\theta}$
and $\chi_{1}, \dots, \chi_{\theta}$; and
\eqref{relations3} holds.

\medbreak
We have found a positive datum  $\mathcal D$ for $\Gamma$
and constructed a  homomorphism of Hopf algebras
$\varphi:  U({\mathcal D}) \to A$. Now
$\gr \varphi: \gr  U({\mathcal D}) \to \gr  A$
is an isomorphism by Theorem \ref{construction};
indeed $\gr \varphi$ is surjective and
the restriction of  $\gr \varphi$ to the first
term of the coradical filtration is injective;
thus $\gr \varphi$ is injective \cite[Th. 5.3.1]{M}.
Hence $\varphi$ is is an isomorphism.
\epf

\begin{obs}
(i). This Theorem can be generalized to the case when $G(A)$
is any abelian group.

(ii). As the proof shows, the condition in (a) that $A$ is a domain
can be replaced by ``$G(A)$ is free abelian of finite rank".

(iii). We believe that the Theorem also holds for generic
infinitesimal braidings.

\end{obs}

\end{document}